January 29, 2010

# EXPONENTIAL POLYNOMIALS, STIRLING NUMBERS, AND EVALUATION OF SOME GAMMA INTEGRALS

Khristo N. Boyadzhiev

Department of Mathematics, Ohio Northern University
Ada, Ohio 45810, USA
k-boyadzhiev@onu.edu

**Abstract**. This article is a short elementary review of the exponential polynomials (also called single-variable Bell polynomials) from the point of view of Analysis. Some new properties are included and several Analysis-related applications are mentioned. At the end of the paper one application is described in details - certain Fourier integrals involving $\Gamma(a+it)$ and $\Gamma(a+it)\Gamma(b-it)$ are evaluated in terms of Stirling numbers.

**1. Introduction.**

We review the exponential polynomials $\phi_n(x)$ and present a list of properties for easy reference. Exponential polynomials in Analysis appear, for instance, in the rule for computing derivatives like $(\frac{d}{dt})^n e^{ae^t}$ and the related Mellin derivatives

$$(x\frac{d}{dx})^n f(x), \quad (\frac{d}{dx}x)^n f(x).$$

Namely, we have

$$(\frac{d}{dt})^n e^{ae^t} = \phi_n(ae^t) e^{ae^t} \tag{1.1}$$

or, after the substitution $x = e^t$,

$$(x\frac{d}{dx})^n e^{ax} = \phi_n(ax) e^{ax}. \tag{1.2}$$

We also include in this review two properties relating exponential polynomials to Bernoulli numbers, $B_k$. One is the semi-orthogonality

$$\int_{-\infty}^{0} \phi_n(x) \phi_m(x) e^{2x} \frac{dx}{x} = (-1)^n \frac{2^{n+m}-1}{n+m} B_{n+m}, \tag{1.3}$$

where the right hand side is zero if $n+m$ is odd. The other property is (2.10).

At the end we give one application. Using exponential polynomials we evaluate the integrals

$$\int_{\mathbb{R}} e^{-it\lambda} t^n \Gamma(a+it) dt, \tag{1.4}$$

and

$$\int_{\mathbb{R}} e^{-it\lambda} t^n \Gamma(a+it) \Gamma(b-it) dt, \tag{1.5}$$

for $n = 0, 1, \ldots$ in terms of Stirling numbers.

## 2. Exponential polynomials

The evaluation of the series

$$S_n = \sum_{k=0}^{\infty} \frac{k^n}{k!}, \quad n = 0, 1, 2, \ldots, \tag{2.1}$$

has a long and interesting history. Clearly, $S_0 = e$, with the agreement that $0^0 = 1$. Several reference books (for instance, [31]) provide the following numbers.

$S_1 = e$, $S_2 = 2e$, $S_3 = 5e$, $S_4 = 15e$, $S_5 = 52e$, $S_6 = 203e$, $S_7 = 877e$, $S_8 = 4140e$.

As noted by H. Gould in [19, p. 93], the problem of evaluating $S_n$ appeared in the Russian journal *Matematicheskii Sbornik*, 3 (1868), p.62, with solution ibid , 4 (1868-9), p. 39.) Evaluations are presented also in two papers by Dobiński and Ligowski. In 1877 G. Dobiński [15] evaluated the first eight series $S_1, \ldots, S_8$ by regrouping:

$$S_1 = \sum_{1}^{\infty} \frac{k}{k!} = 1 + \frac{2}{2!} + \frac{3}{3!} + \ldots = 1 + \frac{1}{1!} + \frac{1}{2!} + \ldots = e$$

$$S_2 = \sum_{1}^{\infty} \frac{k^2}{k!} = 1 + \frac{2^2}{2!} + \frac{3^2}{3!} + \ldots = 1 + \frac{2}{1!} + \frac{3}{2!} + \frac{4}{3!} + \ldots =$$

$$\{1 + \frac{1}{1!} + \frac{1}{2!} + \frac{1}{3!} + \ldots\} + \{\frac{1}{1!} + \frac{2}{2!} + \frac{3}{3!} + \ldots\} = e + S_1 = 2e,$$

and continuing like that to $S_8$. For large $n$ this method is not convenient. However, later that year Ligowski [27] suggested a better method, providing a generating function for the numbers $S_n$

$$e^{e^z} = \sum_{k=0}^{\infty} \frac{e^{kz}}{k!} = \sum_{k=0}^{\infty} \sum_{n=0}^{\infty} \frac{k^n}{k!} \frac{z^n}{n!} = \sum_{n=0}^{\infty} S_n \frac{z^n}{n!}.$$

Further, an effective iteration formula was found

$$S_n = \sum_{j=0}^{n-1} \binom{n-1}{j} S_j$$

by which every $S_n$ can be evaluated starting from $S_1$.

These results were preceded, however, by the work [23] of Johann August Grunert (1797-1872), professor at Greifswalde. Among other things, Grunert obtained formula (2.2) below from which the evaluation of (2.1) follows immediately.

The structure of the series $S_n$ hints at the exponential function. Differentiating the expansion

$$e^x = \sum_{k=0}^{\infty} \frac{x^k}{k!}$$

and multiplying both sides by $x$ we get

$$xe^x = \sum_{k=0}^{\infty} \frac{kx^k}{k!}$$

which, for $x = 1$, gives $S_1 = e$. Repeating the procedure, we find $S_2 = 2e$ from

$$x(xe^x)' = (x + x^2)e^x = \sum_{k=0}^{\infty} \frac{k^2 x^k}{k!}$$

and continuing like that, for every $n = 0, 1, 2, \ldots$, we find the relation

$$(x\frac{d}{dx})^n e^x = \phi_n(x) e^x = \sum_{k=0}^{\infty} \frac{k^n x^k}{k!} \qquad (2.2)$$

where $\phi_n$ are polynomials of degree $n$. Thus,

$$S_n = \phi_n(1) e, \ \forall \ n \geq 0.$$

The polynomials $\phi_n$ deserve a closer look. From the defining relation (2.2) we obtain

$$x(\phi_n e^x)' = x(\phi_n' + \phi_n)e^x = \phi_{n+1} e^x$$

i.e. $\qquad \qquad \qquad \phi_{n+1} = x(\phi_n' + \phi_n) \qquad (2.3)$

which helps to find $\phi_n$ explicitly starting from $\phi_0$,

$$\phi_0(x) = 1$$

$$\phi_1(x) = x$$

$$\phi_2(x) = x^2 + x$$

$$\phi_3(x) = x^3 + 3x^2 + x$$

$$\phi_4(x) = x^4 + 6x^3 + 7x^2 + x$$

$$\phi_5(x) = x^5 + 10x^4 + 25x^3 + 15x^2 + x,$$

and so on. Another interesting relation that easily follows from (2.2) is

$$\phi_{n+1}(x) = x \sum_{k=0}^{n} \binom{n}{k} \phi_k(x) . \tag{2.4}$$

Here is a short proof. Starting from

$$\phi_k(x) e^x = \sum_{j=0}^{\infty} \frac{j^k x^j}{j!}$$

we compute

$$\sum_{k=0}^{n} \binom{n}{k} \phi_k(x) e^x = \sum_{j=0}^{\infty} \frac{x^j}{j!} \sum_{k=0}^{n} \binom{n}{k} j^k = \sum_{j=0}^{\infty} \frac{x^j}{j!} (j+1)^n$$

$$= \frac{1}{x} \sum_{j=0}^{\infty} \frac{(j+1)^{n+1} x^{j+1}}{(j+1)!} = \frac{1}{x} \phi_{n+1}(x) e^x$$

and (2.4) is ready.

From (2.3) and (2.4) one finds immediately

$$\phi_n'(x) = \sum_{k=0}^{n-1} \binom{n}{k} \phi_k(x) \tag{2.5}$$

Obviously, $x = 0$ is a zero for all $\phi_n$, $n > 0$. It can be seen that all the zeros of $\phi_n$ are real, simple, and nonpositive. The nice and short induction argument belongs to Harper [24].

The assertion is true for $n = 1$. Suppose that for some $n$ the polynomial $\phi_n$ has $n$ distinct real non-positive zeros (including $x = 0$). Then the same is true for the function

$$f_n(x) = \phi_n(x) e^x.$$

Moreover, $f_n$ is zero at $-\infty$ and by Rolle's theorem its derivative

$$\frac{d}{dx} f_n = \frac{d}{dx}(\phi_n(x) e^x)$$

has $n$ distinct real negative zeros. It follows that the function

$$\phi_{n+1}(x) e^x = x \frac{d}{dx}(\phi_n(x) e^x)$$

has $n+1$ distinct real non-positive zeros (adding here $x = 0$).

The polynomials $\phi_n$ can be defined also by the exponential generating function (extending Ligowski's formula)

$$e^{x(e^z - 1)} = \sum_{n=0}^{\infty} \phi_n(x) \frac{z^n}{n!} . \qquad (2.6)$$

It is not obvious, however, that the polynomials defined by (2.2) and (2.6) are the same, so we need the following simple statement.

**Proposition 1**. *The polynomials $\phi_n(x)$ defined by (2.2) are exactly the partial derivatives $(\partial/\partial z)^n e^{x(e^z - 1)}$ evaluated at $z = 0$.*

(2.6) follows from (2.2) after expanding the exponential $e^{xe^z}$ in double series and changing the order of summation. A different proof will be given later.

Setting $z = 2k\pi i$, $k = \pm 1, \pm 2, \ldots,$ in the generating function (2.6) one finds

$$e^{2k\pi i} = 1, \; e^{x(e^z - 1)} = e^0 = 1,$$

which shows that the exponential polynomials are linearly dependent

$$1 = \sum_{n=0}^{\infty} \phi_n(x) \frac{(2k\pi i)^n}{n!} \quad \text{or} \quad 0 = \sum_{n=1}^{\infty} \phi_n(x) \frac{(2k\pi i)^n}{n!}, \; k = \pm 1, \pm 2, \ldots, \qquad (2.7)$$

In particular, $\phi_n$ are not orthogonal for any scalar product on polynomials. (However, they have the semi-orthogonality property mentioned in the introduction and proved in Section 4.)

Comparing coefficient for $z$ in the equation

$$e^{(x+y)e^z} = e^{xe^z} e^{ye^z}$$

yields the binomial identity

$$\phi_n(x+y) = \sum_{k=0}^{n} \binom{n}{k} \phi_k(x) \phi_{n-k}(y). \tag{2.8}$$

With $y = -x$ this implies the interesting "orthogonality" relation for $n \geq 1$

$$\sum_{k=0}^{n} \binom{n}{k} \phi_k(x) \phi_{n-k}(-x) = 0. \tag{2.9}$$

Next, let $B_n, n = 0, 1, \ldots,$ be the Bernoulli numbers. Then for $p = 0, 1, \ldots,$ we have

$$\int_0^x \phi_p(t) dt = \frac{1}{p+1} \sum_{k=1}^{p+1} \binom{p+1}{k} B_{p+1-k} \phi_k(x). \tag{2.10}$$

For proof see Example 4 in [3, p.51], or [6].

**Some historical notes**

As already mentioned, formula (2.2) appears in the work of Grunert [23], on p. 260, where he gives also the representation (3.3) below and computes explicitly the first six exponential polynomials. The polynomials $\phi_n$ were studied more systematically (and independently) by S. Ramanujan in his unpublished notebooks. Ramanujan's work is presented and discussed by Bruce Berndt in [3, Part 1, Chapter 3]. Ramanujan, for example, obtained (2.6) from (2.2) and also proved (2.4), (2.5) and (2.10). Later, these polynomials were studied by E.T. Bell [1] and Jacques Touchard [39], [40]. Both Bell and Touchard called them "exponential" polynomials, because of their relation to the exponential function, e.g. (1.1), (1.2), (2.2) and (2.6). This name was used also by Gian-Carlo Rota [34]. As a matter of fact, Bell introduced in [1] a more general class of polynomials of many variables, $Y_{n,k}$, including $\phi_n$ as a particular case. For this reason $\phi_n$ are known also as the single-variable Bell polynomials [13], [20], [21], [41]. These polynomials are also a special case of the actuarial polynomials introduced by Toscano [38] which, on their part, belong to the more general class of Sheffer polynomials [7].

The exponential polynomials appear in a number of papers and in different applications - see [4], [5], [6], [29], [32], [33], [34] and the references therein. In [35] they appear on p. 524 as the horizontal generating functions of the Stirling numbers of the second kind (see below (3.3)).

The numbers

$$b_n = \phi_n(1) = \frac{1}{e} S_n \qquad (2.11)$$

are sometimes called exponential numbers, but a more established name is Bell numbers. They have interesting combinatorial and analytical applications [2], [8], [14], [16], [20], [21], [25], [30], [37], [38]. An extensive list of 202 references for Bell numbers is given in [18].

We note that equation (2.2) can be used to extend $\phi_n$ to $\phi_z$ for any complex number $z$ by the formula

$$\phi_z(x) = e^{-x} \sum_{k=0}^{\infty} \frac{k^z x^k}{k!} \qquad (2.12)$$

(Butzer et al. [9], [10]). The function appearing here is an interesting entire function in both variables, $x$ and $z$. Another possibility is to study the polyexponential function

$$e_s(x, \lambda) = \sum_{n=0}^{\infty} \frac{x^n}{n!(n+\lambda)^s}, \qquad (2.13)$$

where $\text{Re}\,\lambda > 0$. When $s$ is a negative integer, the polyexponential can be written as a finite linear combination of exponential polynomials (see [6]).

**3. Stirling numbers and Mellin derivatives**

The iteration formula (2.3) shows that all polynomials $\phi_n$ have positive integer coefficients. These coefficients are the Stirling numbers of the second kind $\{{n \atop k}\}$ (or $S(n,k)$) - see [12], [14], [17], [22], [26], [35]. Given a set of $n$ elements, $\{{n \atop k}\}$ represents the number of ways this set can be partitioned into $k$ nonempty subsets ($0 \leq k \leq n$). Obviously, $\{^n_1\} = 1$, $\{^n_n\} = 1$ and a short computation gives $\{^n_2\} = 2^{n-1} - 1$. For symmetry one sets $\{^0_0\} = 1$, $\{^n_0\} = 0$. The definition of $\{^n_k\}$ implies the property

$$\{^{n+1}_k\} = k\{^n_k\} + \{^n_{k-1}\} \tag{3.1}$$

(see p.259 in [22]) which helps to compute all $\{^n_k\}$ by iteration. For instance,

$$\{^n_3\} = (3^{n-1} - 2^n + 1)/2 .$$

A general formula for the Stirling numbers of the second kind is

$$\{^n_k\} = \frac{1}{k!} \sum_{j=1}^{k} (-1)^{k-j} \binom{k}{j} j^n . \tag{3.2}$$

**Proposition 2.** *For every* $n = 0, 1, 2, \ldots$

$$\phi_n(x) = \{^n_0\} + \{^n_1\}x + \{^n_2\}x^2 + \ldots + \{^n_n\}x^n = \sum_{k=0}^{n} \{^n_k\}x^k . \tag{3.3}$$

The proof is by induction and is left to the reader. Setting here $x = 1$ we come to the well-known representation for the numbers $S_n$

$$S_n = e(\{^n_0\} + \{^n_1\} + \{^n_2\} + \ldots + \{^n_n\}).$$

It is interesting that formula (3.3) is very old - it was obtained by Grunert [23, p 260] together with the representation (3.2) for the coefficients which are called now Stirling numbers of the second kind. In fact, coefficients of the form

$$\{^n_k\}k!$$

appear in the computations of Euler - see [17].

It is good to note that the representation (3.2) quickly follows from (3.3) and (2.2). First we write

$$\sum_{k=0}^{n} \{^n_k\}x^k = e^{-x} \sum_{k=0}^{\infty} \frac{k^n x^k}{k!} = \{\sum_{j=0}^{\infty} \frac{(-1)^j x^j}{j!}\}\{\sum_{k=0}^{\infty} \frac{k^n x^k}{k!}\},$$

then we multiply the two series by Cauchy's rule and compare coefficients. Thus we come to (3.2). This proof shows very well that the right-hand side in (3.2) is zero when $k > n$.

Next we turn to some special differentiation formulas. Let $D = d/dx$.

**Mellin derivatives**. It is easy to see that the first equality in (2.2) extends to equation (1.2), where $a$ is an arbitrary complex number i.e.,

$$(xD)^n e^{ax} = \phi_n(ax) e^{ax}$$

by the substitution $x \to ax$. Even further, this extends to

$$(xD)^n e^{ax^p} = p^n \phi_n(ax^p) e^{ax^p} \qquad (3.4)$$

for any $a, p$ and $n = 0, 1, \ldots$ (simple induction and (2.3)). Again by induction, it is easy to prove that

$$(xD)^n f(x) = \sum_{k=0}^{n} \{{}_k^n\} x^k D^k f(x). \qquad (3.5)$$

for any $n$-times differentiable function $f$. This formula was obtained by Grunert [23, pp 257-258] (see also p. 89 in [19], where a proof by induction is given).

As we know the action of $xD$ on exponentials, formula (3.5) can be "discovered" by using Fourier transform. Let

$$F[f](t) = \int_{\mathbb{R}} e^{-ixt} f(x) \, dx \qquad (3.6)$$

be the Fourier transform of some function $f$. Then

$$f(x) = \frac{1}{2\pi} \int_{\mathbb{R}} e^{ixt} F[f](t) \, dx,$$

$$(xD)^n f(x) = \frac{1}{2\pi} \int_{\mathbb{R}} e^{ixt} \phi_n(ixt) F[f](t) \, dx = \qquad (3.7)$$

$$\sum_{k=0}^{n} \{{}_k^n\} x^k F^{-1}[(it)^k F[f]](x) = \sum_{k=0}^{n} \{{}_k^n\} x^k D^k f(x) \quad .$$

Next we turn to formula (1.1) and explain its relation to (1.2). If we set $x = e^t$, then for any differentiable function $f$

$$\frac{d}{dt}f = (\frac{d}{dx}f)\frac{dx}{dt} = (\frac{d}{dx}f)e^t = (xD)f$$

and we see that (1.1) and (1.2) are equivalent.

$$(\frac{d}{dt})^n e^{ae^t} = (xD)^n e^{ax} = \phi_n(ax)e^{ax} = \phi_n(ae^t)e^{ae^t} \tag{3.8}$$

*Proof* of Proposition 1. We apply (1.1) to the function $f_x(z) = e^{x(e^z-1)} = e^{xe^z}e^{-x}$

$$(\frac{d}{dz})^n f_x(z) = \phi_n(xe^z)f_x(z)$$

From here, with $z = 0$

$$(\frac{d}{dz})^n f_x(z)|_{z=0} = \phi_n(x)$$

as needed.

Now we list some simple operational formulas. Starting from the obvious relation

$$(xD)^n x^k = k^n x^k, \; n = 0, 1, ..., k \in \mathbb{R}, \tag{3.9}$$

for any function of the form

$$f(t) = \sum_{n=0}^{\infty} a_n t^n, \tag{3.10}$$

we define the differential operator

$$f(xD) = \sum_{n=0}^{\infty} a_n (xD)^n$$

with action on functions $g(x)$,

$$f(xD)g(x) = \sum_{n=0}^{\infty} a_n (xD)^n g(x). \tag{3.11}$$

When $g(x) = x^k$, (3.9) and (3.11) show that

$$f(xD)x^k = \sum_{n=0}^{\infty} a_n k^n x^k = f(k) x^k.$$

If now
$$g(x) = \sum_{k=0}^{\infty} c_k x^k$$

is a function analytical in a neighborhood of zero, the action of $f(xD)$ on this function is given by

$$f(xD)g(x) = \sum_{k=0}^{\infty} c_k f(k) x^k, \tag{3.12}$$

provided the series on the right side converges. When $f$ is a polynomial, formula (3.12) helps to evaluate series like

$$\sum_{k=0}^{\infty} c_k f(k) x^k$$

in a closed form. This idea was exploited by Schwatt [36] and more recently by the present author in [4]. For instance, when $g(x) = e^x$ equation (3.12) becomes

$$\sum_{k=0}^{\infty} f(k) \frac{x^k}{k!} = e^x \sum_{n=0}^{\infty} a_n \phi_n(x). \tag{3.13}$$

As shown in [4] this series transformation can be used for asymptotic series expansions of certain functions.

**Leibniz Rule.** The higher order Mellin derivative $(xD)^n$ satisfies the Leibniz rule

$$(xD)^n(fg) = \sum_{k=0}^{n} \binom{n}{k} [(xD)^{n-k}f][(xD)^k g] \tag{3.14}$$

The proof is easy, by induction, and is left to the reader. We shall use this rule to prove the following proposition.

**Proposition 3.** For all $n, m = 0, 1, 2, \ldots$

$$\phi_{n+m}(x) = \sum_{k=0}^{n} \sum_{j=0}^{m} \binom{n}{k} \{{}^m_j\} j^{n-k} x^j \phi_k(x) \tag{3.15}$$

*Proof.*

$$\phi_{n+m}(x) = (xD)^{n+m} e^x = (xD)^n (xD)^m e^x = (xD)^n (\phi_m(x) e^x),$$

which by the Leibniz rule (3.14) equals

$$\sum_{k=0}^{n} \binom{n}{k} [(xD)^{n-k} \phi_m(x)] [(xD)^k e^x].$$

Using (3.2) and (3.9) we write

$$(xD)^{n-k} \phi_m(x) = \sum_{j=0}^{m} \{{}^m_j\} j^{n-k} x^j,$$

and since also

$$(xD)^k e^x = \phi_k(x) e^x,$$

we obtain (3.15) from (3.14). The proof is completed.

Setting $x = 1$ in (3.14) yields an identity for the Bell numbers.

$$b_{n+m} = \sum_{k=0}^{n} \sum_{j=0}^{m} \binom{n}{k} \{{}^m_j\} j^{n-k} b_k. \tag{3.16}$$

This identity was recently published by Spivey [37], who gave a combinatorial proof. After that Gould and Quaintance [21] obtained the generalization (3.15) together with two equivalent versions. The proof in [21] is different from the one above.

Using the Leibniz rule for $xD$ we can prove also the following extension of property (2.9)

**Proposition 4**. For any two integers $n, m \geq 0$

$$(xD)^n \phi_m(x) = \sum_{k=0}^{m} \{{}^m_k\} k^n x^k = \sum_{k=0}^{n} \binom{n}{k} \phi_{m+k}(x) \phi_{n-k}(-x). \tag{3.17}$$

The proof is simple. Just compute

$$(xD)^n \phi_m(x) = (xD)^n[(e^{-x})(\phi_m(x) e^x)]$$

$$= \sum_{k=0}^{n} \binom{n}{k} [(xD)^{n-k} e^{-x}] [(xD)^k (\phi_m(x) e^x)]$$

and (3.17) follows from (1.2).

For completeness we mention also the following three properties involving the operator $Dx$. Proofs and details are left to the reader.

$$(Dx)^n e^{ax} = \frac{\phi_{n+1}(ax)}{ax} e^{ax}, \tag{3.18}$$

$$(Dx)^n f(x) = \sum_{k=0}^{n} \left\{{n+1 \atop k+1}\right\} x^k D^k f(x), \qquad (3.19)$$

and
$$f(Dx) g(x) = \sum_{k=0}^{\infty} c_k f(k+1) x^k, \qquad (3.20)$$

analogous to (1.2). (3.5), and (3.12) correspondingly

For a comprehensive study of the Mellin derivative we refer to [11].

**More Stirling numbers.** The polynomials $\phi_n$, $n = 0, 1, \ldots$, form a basis in the linear space of all polynomials. Formula (3.3) shows how this basis is expressed in terms of the standard basis $1, x, x^2, \ldots, x^n, \ldots$,. We can solve for $x^k$ in the equations (3.3) and express the standard basis in terms of the exponential polynomials

$$1 = \phi_0$$

$$x = \phi_1$$

$$x^2 = -\phi_1 + \phi_2$$

$$x^3 = 2\phi_1 - 3\phi_2 + \phi_3$$

$$x^4 = -6\phi_1 + 11\phi_2 - 6\phi_3 + \phi_4,$$

etc. The coefficients here are also special numbers. If we write

$$x^n = \sum_{k=0}^{n} (-1)^{n-k} \left[{n \atop k}\right] \phi_k \qquad (3.21)$$

then $\left[{n \atop k}\right]$ are the (absolute) Stirling numbers of first kind, as defined in [22]. (The numbers $\left[{n \atop k}\right]$ are non-negative. The symbol $s(n,k) = (-1)^{n-k} \left[{n \atop k}\right]$ is used for Stirling numbers of the first kind with changing sign - see [14], [18] and [26] for more details.) $\left[{n \atop k}\right]$ is the number of ways to arrange $n$ objects into $k$ cycles. According to this interpretation,

$$\left[{n \atop k}\right] = (n-1)\left[{n-1 \atop k}\right] + \left[{n-1 \atop k-1}\right], \quad n \geq 1.$$

**4. Semi-orthogonality of $\phi_n$**

**Proposition 5.** For every $n, m = 1, 2, \ldots$, we have

$$\int_0^\infty \phi_n(-x)\phi_m(-x) e^{-2x} \frac{dx}{x} = (-1)^{n-1} \frac{2^{n+m}-1}{n+m} B_{n+m}. \tag{4.1}$$

Here $B_k$ are the Bernoulli numbers. Note that the right hand side is zero when $k+m$ is odd, as all Bernoulli numbers with odd indices $> 1$ are zeros.

Using the representation (3.3) in (4.1) and integrating termwise we obtain an equivalent form of (4.1)

$$\sum_{k=0}^n \sum_{j=0}^m (-1)^{k+j} \left\{{n \atop k}\right\}\left\{{m \atop j}\right\} \frac{(k+j-1)!}{2^{k+j}} = (-1)^{n-1} \frac{2^{n+m}-1}{n+m} B_{n+m}. \tag{4.2}$$

This (double sum) identity extends the known identity [22, p.317, Problem 6.76]

$$\sum_{j=0}^m (-1)^{j+1} \left\{{m \atop j}\right\} \frac{j!}{2^{j+1}} = \frac{2^{m+1}-1}{m+1} B_{m+1}. \tag{4.3}$$

Namely, (4.3) results from (4.2) for $n = 1$. The presence of $(-1)^{n-1}$ at the right hand side in (4.1) is not a "break of symmetry", because when $n+m$ is even, then $n$ and $m$ are both even or both odd.

*Proof of the proposition.* Starting from

$$\Gamma(z) = \int_0^\infty x^{z-1} e^{-x} dx \tag{4.4}$$

we set $x = e^\lambda$, $z = a + it$, to obtain the representation

$$\Gamma(a+it) = \int_{-\infty}^{+\infty} e^{i\lambda t} e^{a\lambda} e^{-e^\lambda} d\lambda, \tag{4.5}$$

which is a Fourier transform integral. The inverse transform is

$$e^{a\lambda} e^{-e^\lambda} = \frac{1}{2\pi} \int_{\mathbb{R}} e^{-i\lambda t} \Gamma(a+it) dt. \tag{4.6}$$

When $a = 1$ this is

$$-e^\lambda e^{-e^\lambda} = \frac{d}{d\lambda} e^{-e^\lambda} = \frac{-1}{2\pi} \int_{\mathbb{R}} e^{-i\lambda t} \Gamma(1+it) dt. \tag{4.7}$$

Differentiating (4.7) $n-1$ times for $\lambda$ we find

$$(\frac{d}{d\lambda})^n e^{-e^\lambda} = \phi_n(-e^\lambda) e^{-e^\lambda} = \frac{-1}{2\pi} \int_{\mathbb{R}} e^{-i\lambda t} (-it)^{n-1} \Gamma(1+it) dt, \tag{4.8}$$

and Parceval's formula yields the equation

$$\int_{\mathbb{R}} \phi_n(-e^\lambda) \phi_m(-e^\lambda) e^{-2e^\lambda} d\lambda = \frac{1}{2\pi} \int_{\mathbb{R}} (-it)^{n-1} (it)^{m-1} |\Gamma(1+it)|^2 dt$$

or, with $x = e^\lambda$

$$\int_0^\infty \phi_n(-x) \phi_m(-x) e^{-2x} \frac{dx}{x} = \frac{(-1)^n i^{n+m}}{2\pi} \int_{\mathbb{R}} t^{n+m-2} \frac{\pi t}{\sinh(\pi t)} dt. \tag{4.9}$$

The right hand side is 0 when $n + m$ is odd. When $n + m$ is even, we use the integral [31, p.351]

$$\int_0^\infty \frac{t^{2p-1}}{\sinh(\pi t)} dt = \frac{2^{2p}-1}{2p}(-1)^{p=1} B_{2p} \qquad (4.10)$$

to finish the proof.

Property (4.1) resembles the semi-orthogonal property of the Bernoulli polynomials

$$\int_0^1 B_n(x) B_m(x) dx = (-1)^{n-1} \frac{n!\, m!}{(n+m)!} B_{n+m}, \qquad (4.11)$$

- see, for instance, [35, p.530].

## 5. Gamma integrals.

We use the technique in the previous section to compute certain Fourier integrals and evaluate the moments of $\Gamma(a+it)$ and $\Gamma(a+it)\Gamma(b-it)$.

**Proposition 6**. For every $n = 0, 1, \ldots$ and $a, b > 0$ we have

$$\int_\mathbb{R} e^{-i\mu t} t^n\, \Gamma(a+it)\, \Gamma(b-it)\, dt \qquad (5.1)$$

$$= i^n 2\pi\, e^{-b\mu} \sum_{k=0}^n \sum_{m=0}^k \binom{n}{k} \left\{\begin{matrix}k\\m\end{matrix}\right\} (-1)^m a^{n-k} \frac{\Gamma(a+b+m)}{(1+e^{-\mu})^{a+b+m}} ;$$

$$\int_\mathbb{R} e^{-i\lambda t} t^n \Gamma(a+it)\, dt \qquad (5.2)$$

$$= i^n 2\pi e^{a\lambda} e^{-e^\lambda} \sum_{k=0}^{n} \binom{n}{k} a^{n-k} \sum_{m=0}^{k} \left\{{k \atop m}\right\} (-1)^m e^{\lambda m}.$$

In particular, when $\lambda = \mu = 0$, we obtain the moments

$$G_n(a,b) \equiv \int_{\mathbb{R}} t^n \, \Gamma(a+it) \Gamma(b-it) \, dt \qquad (5.3)$$

$$= i^n \pi \sum_{k=0}^{n} \sum_{m=0}^{k} \binom{n}{k} \left\{{k \atop m}\right\} (-1)^m a^{n-k} \frac{\Gamma(a+b+m)}{2^{a+b+m-1}},$$

$$G_n(a) \equiv \int_{\mathbb{R}} t^n \, \Gamma(a+it) \, dt = \frac{2\pi i^n}{e} \sum_{k=0}^{n} \sum_{m=0}^{k} \binom{n}{k} \left\{{k \atop m}\right\} (-1)^m a^{n-k}. \qquad (5.4)$$

When $n = 0$ in (5.1) we have the known integral

$$\int_{\mathbb{R}} e^{-i\mu t} \Gamma(a+it) \, \Gamma(b-it) \, dt = 2\pi \Gamma(a+b) e^{-b\mu} (1+e^{-\mu})^{-a-b}, \qquad (5.5)$$

which can be found in the form of an inverse Mellin transform in [28].

*Proof.* Using again equation (4.6)

$$e^{a\lambda} e^{-e^\lambda} = \frac{1}{2\pi} \int_{\mathbb{R}} e^{-i\lambda t} \Gamma(a+it) \, dt \qquad (5.6)$$

we differentiate both side $n$ times

$$\left(\frac{d}{d\lambda}\right)^n [e^{a\lambda} e^{-e^\lambda}] = \frac{1}{2\pi} \int_{\mathbb{R}} e^{-i\lambda t} (-it)^n \Gamma(a+it) \, dt,$$

and then, according to the Leibniz rule and (1.1) the left hand side becomes.

$$\left(\frac{d}{d\lambda}\right)^n [e^{a\lambda} e^{-e^\lambda}] = e^{a\lambda} e^{-e^\lambda} \sum_{k=0}^{n} \binom{n}{k} \phi_k(-e^\lambda) a^{n-k}.$$

Therefore,

$$e^{a\lambda} e^{-e^\lambda} \sum_{k=0}^{n} \binom{n}{k} \phi_k(-e^\lambda) a^{n-k} = \frac{1}{2\pi} \int_{\mathbb{R}} e^{-i\lambda t} (-it)^n \Gamma(a+it)\, dt, . \qquad (5.7)$$

and (5.2) follows from here.

Replacing $\lambda$ by $\lambda - \mu$ we write (5.6) in the form

$$e^{b\lambda} e^{-b\mu} e^{-e^\lambda e^{-\mu}} = \frac{1}{2\pi} \int_{\mathbb{R}} e^{-i\lambda t} e^{i\mu t} \Gamma(b+it)\, dt, \qquad (5.8)$$

and then Parceval's formula for Fourier integrals applied to (5.7) and (5.8) yields

$$e^{-b\mu} \sum_{k=0}^{n} \binom{n}{k} a^{n-k} \int_{\mathbb{R}} e^{(a+b)\lambda} e^{-e^\lambda(1+e^{-\mu})} \phi_k(-e^\lambda)\, d\lambda \qquad (5.9)$$

$$= \frac{(-i)^n}{2\pi} \int_{\mathbb{R}} e^{-i\mu t} t^n \Gamma(a+it)\, \Gamma(b-it)\, dt.$$

Returning to the variable $x = e^\lambda$ we write this in the form

$$\frac{1}{2\pi} \int_{\mathbb{R}} e^{-i\mu t} t^n\, \Gamma(a+it)\, \Gamma(b-it)\, dt \qquad (5.10)$$

$$= i^n e^{-b\mu} \sum_{k=0}^{n} \binom{n}{k} a^{n-k} \int_0^\infty \phi_k(-x)\, x^{a+b-1} e^{-x(1+e^{-\mu})}\, dx$$

$$= i^n e^{-b\mu} \sum_{k=0}^{n} \sum_{j=0}^{k} \binom{n}{k} \left\{ {k \atop j} \right\} a^{n-k} (-1)^j \int_0^\infty x^{a+b+j-1} e^{-x(1+e^{-\mu})}\, dx$$

$$= i^n e^{-b\mu} \sum_{k=0}^{n} \sum_{j=0}^{k} \binom{n}{k} \left\{ {k \atop j} \right\} a^{n-k} (-1)^j \frac{\Gamma(a+b+j)}{(1+e^{-\mu})^{a+b+j}}$$

which is (5.1). The proof is complete.

Next, we observe that for any polynomial

$$p(t) = \sum_{n=0}^{m} a_n t^n \qquad (5.11)$$

one can use (5.4) to write the following evaluation

$$\int_{\mathbb{R}} p(t) \, \Gamma(a+it) \, dt = \sum_{n=0}^{m} a_n G_n(a). \qquad (5.12)$$

In particular, when $a = 1$ we have

$$G_n(1) = 2\pi i^n e^{-1} \phi_{n+1}(-1), \qquad (5.13)$$

and therefore,

$$\int_{\mathbb{R}} p(t) \, \Gamma(1+it) \, dt = \frac{2\pi}{e} \sum_{n=0}^{m} a_n i^n \phi_{n+1}(-1). \qquad (5.14)$$

More applications can be found in the recent papers [4], [5] and [6].